\documentclass[11pt]{amsart}
\usepackage{amssymb}

\usepackage{amscd}
\usepackage{amsthm}
\usepackage[dvips]{graphicx}
\usepackage[all]{xy}

\newtheorem{Thm}{Theorem}[section]
\newtheorem{Lem}[Thm]{Lemma}
\newtheorem{Def}[Thm]{Definition}
\newtheorem{Cor}[Thm]{Corollary}
\newtheorem{Prop}[Thm]{Proposition}

\newtheorem{Rem}[Thm]{Remark}

\title[Gerstenhaber Lie bracket]{On the Hochschild cohomology ring of the quaternion group of order eight in characteristic two  }

\author{Alexander Ivanov}

\author{Sergei O. Ivanov}

\author{Yury Volkov }

\author{Guodong Zhou }

\address{ Alexander Ivanov
\newline Department of Higher Algebra and Number Theory
 \newline Faculty of Mathematics and Mechanics
\newline Saint Petersburg State University
\newline  Universitetsky prospekt 28
\newline Saint Petersburg, Peterhof, 198504
 \newline Russia}
\email{a.a.ivanov.spb@gmail.com}

\address{Sergei O. Ivanov
\newline Chebyshev Laboratory, 
\newline St. Petersburg State University, 
\newline 14th Line, 29b, 
\newline Saint Petersburg, 199178 
\newline Russia
}
\email{ivanov.s.o.1986@gmail.com }

\address{ Yury Volkov
\newline Department of Higher Algebra and Number Theory
 \newline Faculty of Mathematics and Mechanics
\newline Saint Petersburg State University
\newline  Universitetsky prospekt 28
\newline Saint Petersburg, Peterhof, 198504
 \newline Russia}
\email{ wolf86\_666@list.ru}

\address{Guodong Zhou
 \newline Department of Mathematics
 \newline Shanghai Key laboratory of PMMP
\newline East China Normal University
\newline  Dong Chuan Road 500
\newline Shanghai 200241
 \newline P.R.China}
  \email{gdzhou@math.ecnu.edu.cn}

\date{version of \today}

\newenvironment{Proof}[1][Proof]{\begin{trivlist}
\item[\hskip \labelsep {\bfseries #1}]}{\flushright
$\Box$\end{trivlist}}

\newcommand{\ot}{\otimes}
\newcommand{\comp}{\mathop{\raisebox{+.3ex}{\hbox{$\scriptstyle\circ$}}}}
\newcommand{\xto}{\xrightarrow}

\newcommand{\lra}{\longrightarrow}

\newcommand{\ra}{\rightarrow}
\newcommand{\sdp}{\times\kern-.2em\vrule height1.1ex depth-.05ex}
\newcommand{\epi}{\lra \kern-.8em\ra}

\newcommand{\ol}{\overline}

 \normalbaselines
\thanks{The first, the second and the third authors are
supported by RFBR grant 13-01-00902\_ a}

\thanks{The first and the second authors are
supported by RF Presidental grant MD-381.2014.1}

\thanks{The second author is  supported by the Chebyshev Laboratory  (Department of Mathematics and Mechanics, St. Petersburg State University)  under RF Government grant 11.G34.31.0026 and by JSC ``Gazprom Neft''.}

\thanks{The fourth  author is
supported by  the exchange program STIC-Asie 'ESCAP' financed by the French Ministry of Foreign Affairs, by Shanghai Pujiang
Program (No. 13PJ1402800),  by National Natural Science Foundation of China (No. 11301186) and by the Doctoral Fund of Youth Scholars of Ministry of Education of China (No. 20130076120001). }

\def\T{{\rm T}}

\begin{document}

\renewcommand{\thefootnote}{\alph{footnote}}
\renewcommand{\thefootnote}{\alph{footnote}}
\setcounter{footnote}{-1}

\renewcommand{\thefootnote}{\alph{footnote}}
\setcounter{footnote}{-1} 

\begin{abstract}  Let $k$ be an algebraically closed field of
characteristic two and let $Q_8$ be the quaternion group of order
$8$. We determine  the Gerstenhaber Lie algebra structure and the
Batalin-Vilkovisky structure on the Hochschild cohomology ring
of the group algebra $kQ_8$.
\end{abstract}

\maketitle

\section*{Introduction}

Let $A$ be an associative algebra over a field $k$. The Hochschild
cohomology $HH^*(A)$ of $A$ has a very rich structure. It  is a graded commutative algebra
via the cup product or the Yoneda product, and  it has a graded  Lie bracket  of degree $-1$ so that it becomes  a graded Lie
algebra; these make ${HH}^*(A)$ a Gerstenhaber
algebra \cite{Gerstenhaber}.

During several decades,  a new structure in Hochschild theory has been
extensively studied in topology and mathematical physics, and recently this was introduced into algebra, the so-called Batalin-Vilkovisky structure. Roughly speaking a Batalin-Vilkovisky (aka. BV) structure is an operator on Hochschild cohomology which squares to zero and which, together with the cup product, can express the Lie bracket.  A BV structure exists only on Hochschild cohomology of certain special classes of algebras.
 T.~Tradler first found  that the Hochschild cohomology algebra of a
finite dimensional symmetric algebra, such as a group algebra of a finite group, is a BV algebra \cite{Tradler}; for later proofs,  see e.g.
 \cite{EuSchedler, Menichi}.

 One of the value of BV structure  is that it gives a method to compute the Gerstenhaber Lie bracket which is usually out of reach in practice.  This paper deals with a concrete example.  Let $k$ be an algebraically closed field of
characteristic two and let $Q_8$ be the quaternion group of order
$8$. In this paper, we compute  explicitly   the Gerstenhaber Lie algebra structure and the
Batalin-Vilkovisky structure on the Hochschild cohomology ring
of the group algebra $kQ_8$.
The Hochschild cohomology ring of $kQ_8$ was calculated by A. I. Generalov in \cite{Generalov} using a minimal projective bimodule  resolution of $kQ_8$. Since the Gerstenhaber Lie bracket is defined using the bar resolution, one needs to find the comparison morphisms between the (normalized) bar resolution and the resolution of Generalov. To this end, we use an effective method employing the notion of weak self-homotopy, recently popularized by J. Le and the fourth author (\cite{LeZhou}).

\section{Hochschild (co)homology}

The cohomology theory of associative algebras was introduced by
Hochschild (\cite{Hochschild}).  The Hochschild cohomology ring of
a $k$-algebra is a Gerstenhaber algebra, which was first
discovered by Gerstenhaber in \cite{Gerstenhaber}. Let us recall
his construction here. Given a $k$-algebra $A$, its Hochschild
cohomology groups are defined as $HH^n(A)\cong
\mathrm{Ext}^n_{A^e}(A, A)$ for $n \geq 0$, where $A^e=A \otimes
A^{\mathrm{op}}$ is the enveloping algebra of $A$. There is a
projective resolution of $A$ as an $A^e$-module

\[
\mathrm{Bar}_*(A) \colon \cdots \to A^{\otimes (r+2)}\xto{d_{r}}
A^{\otimes (r+1)}\to \cdots \to A^{\otimes 3}\xto{d_1} A^{\otimes
2}(\stackrel{d_0=\mu}{\to} A),
\]
where $\mathrm{Bar}_r(A):=A^{\otimes (r+2)}$ for $r\geq 0$, the
map $\mu: A\otimes A\to A$ is the multiplication of $A$, and $d_r$
is defined by

\[
d_r(a_0\otimes a_1\otimes\cdots\otimes a_{r+1})=\sum_{i=0}^r(-1)^i
a_0 \otimes \cdots\otimes a_{i-1}\otimes a_ia_{i+1}\otimes
a_{i+2}\otimes \cdots\otimes a_{r+1}
\]
for all $a_0, \cdots, a_{r+1} \in A$. This is usually called the
(unnormalized) bar resolution of $A$. The normalized version
$\overline{\mathrm{Bar}}_*(A)$ is given by
$\overline{\mathrm{Bar}}_r(A)=A\otimes \overline{A}^{\otimes
r}\otimes A$,  where $\overline{A}=A/(k\cdot 1_A)$, and with the
induced differential from that of $\mathrm{Bar}_*(A)$.

The complex which is used to compute the Hochschild cohomology is
$C^*(A)=\mathrm{Hom}_{A^e}(\mathrm{Bar}_*(A), A)$. Note that for
each $r\geq 0$, $C^r(A)=\mathrm{Hom}_{A^e}(A^{\otimes (r+2)},
A)\cong\mathrm{Hom}_k(A^{\otimes r}, A)$. If $f\in C^r(A)$, then the expression $f(a)$ makes sense for $a\in A^{\otimes (r+2)}$ and $a\in A^{\otimes r}$ simultaneously.
 We identify $C^0(A)$ with $A$.  Thus $C^*(A)$ has the following form:
$$C^*(A) \colon A \xto{\delta^0} \mathrm{Hom}_k(A, A) \to \cdots \to \mathrm{Hom}_k(A^{\otimes r}, A)
\xto{\delta^r} \mathrm{Hom}_k(A^{\otimes (r+1)}, A)\to \cdots .$$
Given $f$ in $\mathrm{Hom}_k(A^{\otimes r}, A)$, the map
$\delta^r(f)$ is defined by sending $a_1\otimes \cdots\otimes
a_{r+1}$ to
\[
a_1 \cdot f(a_2\otimes \cdots\otimes
a_{r+1})+\sum_{i=1}^r(-1)^if(a_1\otimes \cdots \otimes
a_{i-1}\otimes a_ia_{i+1} \otimes a_{i+2}\otimes \cdots\otimes
a_{r+1})+(-1)^{r+1}f(a_1\otimes \cdots \otimes a_r)\cdot a_{r+1}.
\]
There is also a normalized version
$\overline{C}^*(A)=\mathrm{Hom}_{A^e}(\overline{\mathrm{Bar}}_*(A),
A)\cong \mathrm{Hom}_k(\overline{A}^{\otimes *}, A)$.

The cup product $\alpha\smile\beta \in
C^{n+m}(A)=\mathrm{Hom}_k(A^{\otimes (n+m)}, A)$ for $\alpha\in
C^n(A)$  and $\beta\in C^m(A)$ is given by
\[
(\alpha\smile\beta)(a_1\otimes \cdots\otimes
a_{n+m}):=\alpha(a_1\otimes\cdots\otimes a_n)\cdot
\beta(a_{n+1}\otimes\cdots\otimes a_{n+m}).
\]
This cup product induces a well-defined product in Hochschild
cohomology
\[
\smile \colon HH^n(A) \times HH^m(A) \longrightarrow HH^{n+m}(A)
\]
which turns the graded $k$-vector space $HH^*(A)=\bigoplus_{n\geq
0}HH^n(A)$ into a graded commutative algebra (\cite[Corollary
1]{Gerstenhaber}).

The Lie bracket is defined as follows. Let $\alpha \in C^n(A)$ and
$\beta \in C^m(A)$. If $n, m\geq 1$, then for $1\leq i\leq n$, set
$\alpha\circ_i \beta \in C^{n+m-1}(A)$ by
\[
(\alpha\circ_i \beta)(a_1\otimes \cdots\otimes
a_{n+m-1}):=\alpha(a_1\otimes \cdots \otimes a_{i-1}\otimes
\beta(a_i\otimes \cdots \otimes a_{i+m-1})\otimes a_{i+m}\otimes
\cdots \otimes a_{n+m-1});
\]
if $ n\geq 1$ and $m=0$, then $\beta\in A$ and for $1\leq i\leq
n$, set
\[
(\alpha\circ_i\beta)(a_1\otimes \cdots\otimes
a_{n-1}):=\alpha(a_1\otimes \cdots \otimes a_{i-1}\otimes \beta
\otimes a_{i }\otimes \cdots \otimes a_{n-1});
\]
for any  other case, set $\alpha\circ_i \beta$ to be zero. Now
define
\[
\alpha\circ \beta
:=\sum_{i=1}^n(-1)^{(m-1)(i-1)}\alpha\circ_i\beta
\]
and
\[
[\alpha,\, \beta] :=\alpha\circ \beta-(-1)^{(n-1)(m-1)}\beta\circ
\alpha.
\]
Note that $[\alpha,\, \beta]\in C^{n+m-1}(A)$. The above $[\
\,,\,\ ]$ induces a well-defined Lie bracket in Hochschild
cohomology
\[
[\ \,,\,\ ] \colon HH^n(A) \times HH^m(A) \longrightarrow
HH^{n+m-1}(A)
\]
such that $(HH^*(A),\, \smile,\, [\ \,,\,\ ])$ is a Gerstenhaber
algebra (\cite{Gerstenhaber}).

The complex used to compute the Hochschild homology $HH_*(A)$ is
$C_*(A)=A\otimes_{A^e} \mathrm{Bar}_*(A)$. Notice that
$C_r(A)=A\otimes_{A^e}A^{\otimes (r+2)}\simeq A^{\otimes (r+1)}$
and the differential $\partial_r: C_r(A)= A^{\otimes (r+1)}\to
C_{r-1}(A)=A^{\otimes r}$  sends $a_0\otimes \cdots \otimes a_{r}$
to $\sum_{i=0}^{r-1} (-1)^i a_0 \otimes \cdots\otimes
a_{i-1}\otimes a_ia_{i+1}\otimes a_{i+2} \otimes \cdots\otimes
a_r+(-1)^r a_r a_0\otimes a_1\otimes \cdots\otimes a_{r-1}.$

There is a  Connes' $\mathfrak{B}$-operator in the Hochschild
homology theory which is defined as follows. For $a_0\otimes
\cdots \otimes a_r\in C_r(A)$, let $\mathfrak{B}(a_0\otimes \cdots
\otimes a_r) \in C_{r+1}(A)$ be
\[
\sum_{i=0}^r (-1)^{ir} 1\otimes a_i\otimes\cdots \otimes
a_r\otimes a_0\otimes \cdots \otimes a_{i-1}+ \sum_{i=0}^r
(-1)^{ir}a_i\otimes 1 \otimes  a_{i+1}\otimes \cdots \otimes
a_r\otimes a_0\otimes \cdots \otimes a_{i-1}.
\]
It is easy to check that $\mathfrak{B}$ is a chain map satisfying
$\mathfrak{B} \comp \mathfrak{B}=0$, which induces an operator
$\mathfrak{B}: HH_r(A)\rightarrow HH_{r+1}(A)$.

All the above constructions, the cup product, the Lie bracket, the
Connes' $\mathfrak{B}$-operator, carry over to normalized
complexes.

\begin{Def}\label{Def-BV-algebra} A {\em Batalin--Vilkovisky  algebra} (BV algebra for short) is a
Gerstenhaber algebra $(A^\bullet,\, \smile,\, [\ \,,\,\ ])$
together with an operator $\Delta\colon A^\bullet \rightarrow
A^{\bullet-1}$ of degree $-1$ such that $\Delta\circ  \Delta=0$
and
\[
[a,\, b]=-(-1)^{(|a|-1)|b|}(\Delta(a\smile b)- \Delta(a)\smile
b-(-1)^{|a|}a\smile \Delta(b))
\]
for   homogeneous  elements $a, b\in A^\bullet$.
\end{Def}

Tradler noticed that the Hochschild cohomology algebra of a
symmetric algebra is a BV algebra \cite{Tradler}, see
also~\cite{Menichi,EuSchedler}. For a symmetric  algebra $A$, he
showed that the $\Delta$-operator on the Hochschild cohomology
corresponds to the Connes' $\mathfrak{B}$-operator on the
Hochschild homology via the duality between the Hochschild
cohomology and the Hochschild homology.

Recall that a finite dimensional $k$-algebra $A$ is called
symmetric if $A$ is isomorphic to its dual $DA=\mathrm{Hom}_k(A,
k)$ as $A^e$-module, or equivalently, if there exists a symmetric
associative non-degenerate bilinear form $\langle\ ,\
\rangle\colon A\times A \rightarrow k$.  This bilinear form
induces a duality between the Hochschild cohomology and the
homology. In fact,
\[
\begin{array}{rcl}\mathrm{Hom}_k(C_*(A),\,
k)&=&\mathrm{Hom}_k(A\otimes_{A^e}\mathrm{Bar}_*(A),\, k)\\
&\cong&
\mathrm{Hom}_{A^e}(\mathrm{Bar}_*(A),\, \mathrm{Hom}_k(A, k))\\
&\cong&
 \mathrm{Hom}_{A^e}(\mathrm{Bar}_*(A),\, A)=C^*(A).\end{array}
\]
Via this duality,  for $n\geq 1$  we obtain     an operator
$\Delta\colon HH^{n}(A)\rightarrow HH^{n-1}(A)$ which is the dual
of Connes' operator.

We recall the following theorem by Tradler.

\begin{Thm} \cite[Theorem 1]{Tradler} \label{Tradler}
With the notation above,  together  with the cup product, the Lie
bracket and the $\Delta$-operator defined above, the Hochschild
cohomology of $A$ is a BV algebra. More precisely, for $\alpha\in
C^n(A)=\mathrm{Hom}_k(A^{\otimes n}, A)$, $\Delta(\alpha) \in
C^{n-1}(A)=\mathrm{Hom}_k(A^{\otimes(n-1)}, A)$ is given by the
equation
\[
\langle \Delta(\alpha)(a_1\otimes\cdots\otimes a_{n-1}),\;
a_n\rangle=\sum_{i=1}^n  (-1)^{i(n-1)} \langle \alpha(a_i\otimes
\cdots\otimes a_{n-1}\otimes a_n\otimes a_1\otimes \cdots\otimes
a_{i-1}), \; 1\rangle
 \]
for $a_1, \cdots, a_n\in A$. The same formula holds also for the
normalized complex $\overline{C}^*(A)$.
\end{Thm}

\section{Constructing comparison
morphisms}\label{Method}

 Let $k$ be a field and let  $B$ be a $k$-algebra.  Given two
left $B$-modules $M$ and $N$, let $P_*$ (resp. $Q_*$) be a projective
resolutions of $M$ (resp. $N$). Then given a homomorphism  of
$B$-modules $f:M\to N$,  it is well known that there exists a
chain map $f_*: P_*\to Q_*$ lifting $f$ (and different lifts are
equivalent up to homotopy). However, sometimes in practice we need
the actual construction of this chain map, called comparison
morphism,  to perform actual computations. This section presents a
method to construct them. The method is not new and it is explained in the book
of Mac Lane; see \cite[Chapter IX Theorem 6.2]{MacLane}.

Our setup is the following.  Suppose that
$$\cdots \longrightarrow P_n\stackrel{d_n^P}{\longrightarrow} P_{n-1}\stackrel{d_{n-1}^P}{\longrightarrow}
\cdots \stackrel{d_1^P}{\longrightarrow} P_0
(\stackrel{d_0^P}{\longrightarrow} M\rightarrow 0)$$ is a projective
resolution of $M$. Then for each $n\geq 0$ there are sets $\{e_{n,i}\}_{i\in X_n}\subset P_n$ and $\{f_{n,i}\}_{i\in X_n}\subset \mathrm{Hom}_B(P_n,B)$ such that $x=\sum_{i\in X_n}f_{n,i}(x)e_{n,i}$ for all $x\in P_n$. Suppose that the second projective resolution
$$\cdots \longrightarrow Q_n\stackrel{d_n^Q}{\longrightarrow} Q_{n-1}\stackrel{d_{n-1}^Q}{\longrightarrow}
\cdots \stackrel{d_1^Q}{\longrightarrow} Q_0
(\stackrel{d_0^Q}{\longrightarrow} N\rightarrow 0)$$
 has a weak self-homotopy in the sense of the following
 definition.

 \begin{Def}\cite{BianZhangZhang} Let $$\cdots Q_n\stackrel{d_n^Q}{\to} Q_{n-1}\stackrel{d_{n-1}^Q}{\to} \cdots
 \stackrel{d_2^Q}{\to}
 Q _1 \stackrel{d_1^Q}{\to} Q_0 \stackrel{d_0^Q}{\to} N\to 0$$ be
 a
  complex. A weak self-homotopy of this complex is a
 collection of $k$-linear maps
 $t_n:
Q_n\rightarrow Q_{n+1}$ for each $n\geq 0$ and $t_{-1}:
M\rightarrow Q_0$ such that for $n\geq 0$,
$t_{n-1}d_n^Q+d_{n+1}^Qt_n=Id_{Q_n}$ and $d_0^Qt_{-1}=Id_N$.

\end{Def}

Now we construct a chain map $f_n: P_n\rightarrow Q_n$ for $n\geq
0$ lifting $f_{-1}=f$. We need to specify the
value of $f_n$ on the elements $e_{n,i}$ for all $i\in X_n$.

For $n=0$, define $f_0(e_{0,i})=t_{-1}fd_0^P(e_{0,i})$. Then
$d_0^Qf_0(e_{0,i})=d_0^Qt_{-1}fd_0^P(e_{0,i})=fd_0^P(e_{0,i})$.

Suppose that we have constructed $f_0, \cdots, f_{n-1}$ such that
for $0\leq i\leq n-1 $,  $d_i^Qf_i=f_{i-1}d_i^P$.  Define $f_n(e_{n,i})=t_{n-1}f_{n-1}d_n^P(e_{n,i})$. It is easy to check that
$$\begin{array}{rcl} d_n^Qf_n(e_{n,i})&=&f_{n-1}d_n^P(e_{n,i}).\end{array}$$ This proves the following
\begin{Prop}
The maps $f_*$ constructed above form a chain map from $P_*$ to
$Q_*$ lifting $f:M\to N$.

\end{Prop}

This result reduces the computation of comparison morphisms to the construction of weak self-homotopies.
It is easy to see that the complex $Q_*$ is exact if and only if there exists a weak self-homotopy of it.   In fact, we can obtain more.
Denote $Z_n=\mathrm{Ker}(d_n)$ for $n\geq 0$ and $Z_{-1}=N$. As vector spaces, one can fix a decomposition of $Q_n=Z_n\oplus Z_{n-1}$ for $n\geq 0$.  Under these identifications, the differential $d_n$ is equal to $\left(\begin{array}{cc} 0 & Id \\ 0 & 0\end{array}\right):  Z_n  \oplus   Z_{n-1}  \to Z_{n-1}  \oplus   Z_{n-2}$ and we can define $t_{-1}=\left(\begin{array}{c} 0 \\ Id \end{array}\right): Z_{-1}\to Z_0\oplus Z_{-1}$ and for $n\geq 0$, $t_n: Z_n\oplus Z_{n-1}\to Z_{n+1}\oplus Z_n$ to be the map $\left(\begin{array}{cc} 0 & 0 \\ Id & 0\end{array}\right)$. Note that our construction has an additional property:

\begin{Lem}\label{AdditionalProperty} For an exact complex of modules over a $k$-algebra, one can always find a weak self-homotopy $\{t_i, i\geq -1\}$ such that $t_{i+1}t_i=0$ for any $i\geq -1$.
\end{Lem}

We are interested in computing Hochschild cohomology of algebras. Let $A$ be a $k$-algebra.  In order to compute Hochschild (co)homology of $A$, one needs a projective resolution of $A$ as a bimodule.   Since this resolution splits as complexes of one-sided modules, one can even choose a weak self-homotopy which are right module homomorphisms and which satisfies the additional property in Lemma~\ref{AdditionalProperty}.

Now let $P_*$ be an $A^e$-projective resolution of $A$. Denote  now $Q_*=\mathrm{Bar}_*(A)$ (or $Q_*=\overline{\mathrm{Bar}}_*(A)$). Let us consider the construction of comparison morphisms $\Psi_*: Q_*\to P_*$ and $\Phi_*: P_*\to Q_*$.

Suppose now that $Q_*=\overline{\mathrm{Bar}}_*(A)$. In this case $Q_*$ has a weak self-homotopy $s_*$ defined by the formula
$$s_n(a_0\otimes a_1\otimes\dots\otimes a_n\otimes 1)=1\otimes a_0\otimes \dots\otimes a_n\otimes 1.$$
Note that $s_{n+1}s_n=0$ for $n\geq -1$, as we are working with the normalized Bar resolution.
  Suppose that the homomorphism $d_n^P$ is defined by the formula
$$
d_n^P(e_{n,i})=\sum\limits_{j\in X_{n-1}}\sum\limits_{p\in T_{n, i,j}}a_pe_{n-1, j}b_p+\sum\limits_{j\in X_{n-1}}\sum\limits_{q\in T_{n, i,j}'}e_{n-1,j}b_q',
$$
where  $a_p\in J_A$ (here $J_A$ is the Jacobson radical of $A$), $b_p, b_q'\in A$,  $T_{n, i,j}$ and $T_{n, i,j}'$ are certain index sets and $e_{n,i}$ as above.

\begin{Lem}\label{Phi}
If $\Phi_*:P_*\rightarrow Q_*$ is the chain map constructed using $s_*$, then
\begin{equation}
\Phi_n(e_{n, i})=1\otimes\sum\limits_{j\in X_{n-1}}\sum\limits_{p\in T_{n, i,j}}a_p\Phi_{n-1}(e_{n-1,j})b_p.
\end{equation}
\end{Lem}

\begin{Proof} By construction, $\Phi_n(e_{n,i})=s_{n-1}\Phi_{n-1}d_n(e_{n, i})$.  Note that $$\Phi_{n-1}(e_{n-1, j})=s_{n-2}\Phi_{n-2}d_{n-1}^P(e_{n-1, j})$$ and thus
$$s_{n-1}\Phi_{n-1}(e_{n-1,j}b_q')=s_{n-1}s_{n-2}\Phi_{n-2}d_{n-1}^P(e_{n-1, j}) b_q'=0.$$  Therefore,
$$\begin{array}{rcl}\Phi_n(e_{n,i})&=&s_{n-1}\Phi_{n-1}d_n(e_{n, i})\\
&=&s_{n-1}\Phi_{n-1}(\sum\limits_{j\in X_{n-1}}\sum\limits_{p\in T_{n, i,j}}a_pe_{n-1, j}b_p)\\
&=&1\otimes\sum\limits_{j\in X_{n-1}}\sum\limits_{p\in T_{n, i,j}}a_p\Phi_{n-1}(e_{n-1,j})b_p.\end{array}$$
\end{Proof}

Let $\mathcal{B}$ be some $k$-basis of $A$ (or $\overline{A}$ in the case $Q_*=\overline{\mathrm{Bar}}_*(A)$). Then  the set $Y_n=\{1\otimes b_1\otimes\dots\otimes b_n\otimes 1\mid b_1,\dots,b_n\in\mathcal{B}\}$ is a basis for $Q_n$  as free $A^e$-module.
Suppose that we have constructed a weak self-homotopy  $t_*$  of $P_*$  such that  $t_{n+1}t_n=0$  and that $t_n$ is a homomorphism of right $A$-modules for all $n\ge -1$.

\begin{Lem}\label{Psi}
 If $\Psi_*:Q_*\rightarrow P_*$ is the chain map constructed using $t_*$, then
\begin{equation}
\Psi_n(1\otimes a_1\otimes\dots\otimes a_n\otimes 1)=t_{n-1}(a_1\Psi_{n-1}(1\otimes a_2\otimes \dots\otimes a_n\otimes 1))
\end{equation}
 for $n\ge 1$ and $a_i\in A$ $(1\le i\le n)$
\end{Lem}

\begin{Proof}
Denote
$$
y=\sum\limits_{i=1}^{n-1}(-1)^i1\otimes a_1\otimes\dots a_{i-1}\otimes a_ia_{i+1}\otimes a_{i+2}\otimes\dots\otimes a_n\otimes 1+(-1)^n1\otimes a_1\otimes\dots\otimes a_n.
$$
As $t_{n-1}t_{n-2}=0$, we have
$$
\begin{array}{rcl}
&&\Psi_n(1\otimes a_1\otimes\dots\otimes a_n\otimes 1)\\
&=&t_{n-1}\Psi_{n-1}d_n^Q(1\otimes a_1\otimes\dots\otimes a_n\otimes 1)=t_{n-1}\Psi_{n-1}(a_1\otimes\dots\otimes a_n\otimes 1+y)\\
&=&t_{n-1}(a_1\Psi_{n-1}(1\otimes a_2\otimes \dots\otimes a_n\otimes 1))+t_{n-1}t_{n-2}\Psi_{n-2}d_{n-1}^Q(y)\\
&=&t_{n-1}(a_1\Psi_{n-1}(1\otimes a_2\otimes \dots\otimes a_n\otimes 1)).
\end{array}
$$

\end{Proof}

Let $V=\oplus_{i=1}^n k x_i$ be a $k$-vector space with basis $\{x_i, 1\leq i\leq n\}$. Let
$A=T(V)/I=k\langle x_1, \cdots, x_n\rangle /I$ be an algebra given by
generators and relations. Then the minimal projective bimodule
resolution of $A$ begins with
\begin{equation} \label{MinimalResolution}\cdots \to  A\otimes R\ot A\stackrel{d_2}{\to} A\otimes V \ot A\stackrel{d_1}{\to}
A\otimes A\stackrel{d_0}{\to} A \to 0,\end{equation} where
\begin{itemize}

\item $V=\oplus_{i=1}^n k x_i$, $R$ is a $k$-complement of $JI+IJ$
in $I$ (thus $R$ is a set of minimal relations), where $J$ is the ideal of $k\langle x_1, \cdots,
x_n\rangle$ generated by $x_1, \cdots, x_n$;

\item $d_0$ is the multiplication of $A$;

\item $d_1$ is induced by $d_1(1\ot x_i\ot 1)=x_i\ot 1-1\ot x_i$ for $1\leq i\leq n$;

\item $d_2$ is induced by the restriction to $R$ of the bimodule
derivation ${\sf C}: TV\to TV\ot V\ot TV$ sending  a path
$x_{i_1}x_{i_2}\cdots x_{i_r}$ (with $1\leq i_1, \cdots, i_r\leq n$) to $\sum_{j=1}^r x_{i_1}\cdots
x_{i_{j-1}}\ot x_{i_j}\ot x_{i_{j+1}}\cdots x_{i_r}$.

\end{itemize}

We shall construct the first three maps of a weak self-homotopy
of this projective resolution, which are moreover right module
homomorphisms.  Let $\mathcal{B}$ be the basis of $A$ formed by
monomials in $x_1, \cdots, x_n$.

 The first two are easy. We define   $t_{-1}=1\ot 1 $ and $t_0(b\ot 1)={\sf C} (b)$ for $b\in
\mathrm{B}$.

For $t_1: A\ot V\ot A\to A\ot R\ot A$, we first fix a vector space
decomposition $TV/I^2=A\oplus I/I^2$. The space $R$, identified
with $I/(JI+IJ)$, generates  $I/I^2$ considered as
$A$-$A$-bimodule. For $b\in \mathcal{B}$, consider $bx_i\in
TV/I^2$, then we can write $bx_i=\sum_{b'\in \mathcal{B}}
\lambda_{b'}b'+\sum_{j} p_j r_j q_j$ with $r_j\in R$ via the
vector space decomposition $TV/I^2=A\oplus I/I^2$. We define
$$t_1(b\ot x_i\ot 1)=\sum_j p_j \ot r_j \ot q_j.$$

\begin{Prop} The above defined maps $t_{-1}, t_0, t_1$ form the first three maps of a weak self-homotopy of the minimal projective bimodule resolution (\ref{MinimalResolution}).

\end{Prop}

\begin{Proof} We have $d_0t_{-1}(1)=d_0(1\ot 1)=1$ and thus $d_0t_{-1}=Id$.

For $b=x_{i_1}x_{i_2}\cdots x_{i_r}\in \mathcal{B}$,  $t_{-1}d_0(b\ot 1)=t_{-1}(b)=1\ot b$, and
$$\begin{array}{rcl} d_1t_0(b\ot 1)&=&d_1{\sf C}(b)\\ &=&d_1(\sum_{j=1}^r x_{i_1}\cdots
x_{i_{j-1}}\ot x_{i_j}\ot x_{i_{j+1}}\cdots x_{i_r})\\
&=&\sum_{j=1}^r x_{i_1}\cdots
x_{i_{j-1}}  x_{i_j}\ot x_{i_{j+1}}\cdots x_{i_r}-\sum_{j=1}^r x_{i_1}\cdots
x_{i_{j-1}}\ot x_{i_j}  x_{i_{j+1}}\cdots x_{i_r}\\
&=& b\ot 1-1\ot b.\end{array}$$
Therefore,  $(d_1t_0+t_{-1}d_0)(b\ot 1)=b\ot 1-1\ot b+1\ot b=b\ot 1$.

Now for $b\in \mathcal{B}$ and $1\leq i\leq n$,
$t_0d_1(b\ot x_i\ot 1)=t_0(bx_i\ot 1-b\ot x_i)$. Recall that via the decomposition $TV/I^2=A\oplus I/I^2$, $bx_i=\sum_{b'\in \mathcal{B}}
\lambda_{b'}b'+\sum_{j} p_j r_j q_j$, so $$t_0d_1(b\ot x_i\ot 1)=t_0(bx_i\ot 1-b\ot x_i)=\sum_{b'\in \mathcal{B}}\lambda_{b'} {\sf C}(b')-{\sf C}(b)x_i.$$
We have also $$d_2t_1(b\ot x_i\ot 1)=d_2(\sum_j p_j\ot r_j\ot q_j)=\sum_j p_j {\sf C}(r_j) q_j.$$

Recall that the bimodule derivation ${\sf C}: TV\to TV\ot V\ot TV$ composed with the sujection $TV\ot V\ot TV\to A\ot V\ot A$ vanishes on $I^2$ and thus induces a well-defined map ${\sf C}: TV/I^2\to A\ot V\ot A$. Furthermore, ${\sf C}$ restricted  to $I/I^2$ is a homomorphism of $A$-$A$-bimodules. This shows that  $\sum_j p_j {\sf C}(r_j) q_j+\sum_{b'\in \mathcal{B}}\lambda_{b'} {\sf C}(b')={\sf C}(bx_i)$ and since ${\sf C}(bx_i)={\sf C}(b)x_i+b\ot x_i\ot 1$, we obtain that $(t_0d_1+d_2t_1)(b\ot x_i\ot 1)=b\ot x_i\ot 1$.

This completes the proof.

\end{Proof}

\section{Weak self-homotopy for $kQ_8$}

Let $k$ be an algebrically closed field of characteristic two. Let
$Q_8$ be the quaternion group of order $8$. Denote by $A=kQ_8$ its
group algebra.  It is well known that $A$ is isomorphic to the
following bounded quiver algebra $kQ/I$:

 $$\xymatrix{\bullet\ar@(ul,dl)[]_{x} \ar@(ur,dr)[]^{y}}$$\\

with relations $$x^2+yxy, y^2+xyx, x^4, y^4.$$\\

 The structure of $A$ can be visualised as follows:

 $$
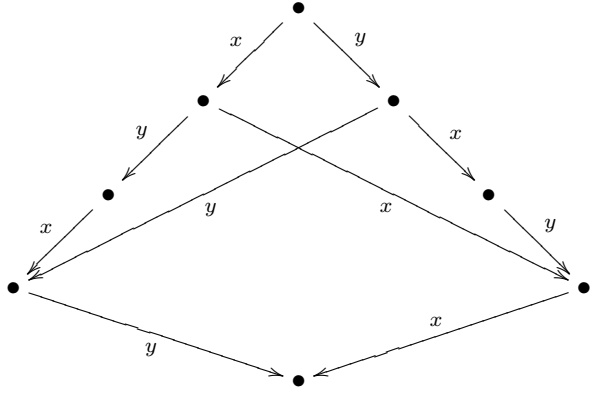
\xymatrix{ & & &\bullet \ar[ld]_x\ar[rd]^y & & & \\
 &&\bullet \ar[ld]_y\ar[rrrrdd]_x && \bullet \ar[lllldd]^y\ar[rd]^x &&\\
 & \bullet \ar[ld]_x& & & & \bullet \ar[rd]^y& \\
 \bullet \ar[rrrd]_y& & & & & &\bullet \ar[llld]_x\\
 & & &\bullet & & &}$$

 A basis of $A$ is given by $\mathcal{B}=\{ 1, x, y, xy, yx, xyx,
 yxy, xyxy\}$. Notice that $\mathcal{B}$ contains a basis of the
 socle of $A$.

  The group algebra  $A$ is a symmetric algebra,
  with respect to the symmetrising form
  $$\langle b_1, b_2\rangle=\left\{\begin{array}{rl} 1 &
  \mathrm{if}\ b_1b_2\in Soc(A)\\
  0 & \mathrm{otherwise}\end{array}\right.$$
  with $b_1, b_2\in \mathcal{B}$.
  The correspondance between elements of $\mathcal{B}$ and its   dual
  basis $\mathcal{B}^*$
  is given by
 $$\begin{array}{ccccccccc} b\in \mathcal{B} &  1 &  x& y&  xy&  yx&  xyx&
 yxy& xyxy\\
 b^*\in \mathcal{B}^* & xyxy & yxy & xyx & xy & yx & y & x & 1\end{array}$$

Since $A$ is an algebra with a DTI-family of relations,  there is
a minimal projective resolution constructed by the second author
in \cite{Sergey}. Let us recall the concrete construction of this
resolution.

After (\cite{Generalov}),   there is an exact sequences of
bimodules  as follows:
$$ (0 \to A\stackrel{\rho}{\to})  A\otimes A \stackrel{d_3}{\to}
A\otimes kQ_1^*\ot A \stackrel{d_2}{\to}  A\otimes kQ_1\ot A
\stackrel{d_1}{\to}A\otimes    A  \stackrel{d_0}{\to} A\to 0$$
where
\begin{itemize}
\item $Q_1=\{x, y\}$ and $Q_1^*=\{r_x, r_y\}$ with $r_x=x^2+yxy$
and $ r_y=y^2+xyx;$

\item  the map $d_0$ is the multiplication of $A$;

\item $d_1(1\ot x\ot 1)=x\ot 1+1\ot x$  and $d_1(1\ot y\ot 1)=y\ot
1+1\ot y$;

\item $d_2(1\ot r_x\ot 1)=1\ot x\ot x+x\ot x\ot 1+ 1\ot y\ot  xy +
  y\ot  x\ot y+     y   x\ot y\ot 1$ \\ and \\
$d_2(1\ot r_y\ot 1)=1\ot y\ot y+y\ot y\ot 1+ 1\ot x\ot  yx +
  x\ot  y\ot x+     xy\ot x\ot 1;$

  \item $d_3(1\ot 1)=x\ot r_x\ot 1 + 1\ot r_x\ot x+ y\ot r_y\ot 1 + 1\ot r_y\ot
  y$;

  \item $\rho(1)=\sum_{b\in \mathcal{B}} b^*\ot b$.

  \end{itemize}

  Using this exact sequence, one can construct a minimal
  projective bimodule resolution of $A$ which is periodic of
  period $4$:

  \begin{itemize}

  \item $P_0=A\otimes A=P_3$, $P_1= A\otimes kQ_1\ot A$ and  $P_2= A\otimes kQ_1^*\ot
  A$;

  \item  $P_4=P_0=A\otimes A$ and $d_4= \rho \circ d_0: P_4\to
  P_3$;

  \item for $n\geq 1$ and
  $i\in \{0, 1, 2, 3\}$, we have
 $P_{4n+i}=P_i$  and $d_{4n+i+1}=d_{i+1}$.

\end{itemize}

We shall establish a weak self-homotopy $\{t_i: P_i\to P_{i+1};
t_{-1}: A\to P_0\}$ over this periodic resolution which are right
module homomorphisms.

 The first two are easy which are
$t_{-1}=1\ot 1 $ and $t_0(b\ot 1)={\sf C} (b)$ for $b\in
\mathrm{B}$, where ${\sf C}: kQ\to kQ\ot kQ_1\ot kQ$ is the
bimodule derivation sending a path $\alpha_1\cdots \alpha_n$ with
$\alpha_1, \cdots, \alpha_n\in Q_1$ to $\sum_{i=1}^n
\alpha_1\cdots \alpha_{i-1}\ot \alpha_i\ot \alpha_{i+1}\cdots
\alpha_n$.

The map $t_1: P_1\to P_2$ is given by
$$\begin{array}{lcl} t_1(1\ot x\ot 1)&=&0, \\

  t_1(x\ot x\ot 1)&=&1\ot r_x\ot
1,\\
t_1(y\ot x\ot 1)&=&0,\\  t_1(xy\ot x\ot 1)&=&0, \\

 t_1(yx\ot x\ot
1)&=&y\ot r_x\ot 1  +xy\ot r_x\ot y +1\ot r_y\ot xy,\\

 t_1(xyx\ot
x\ot 1)&=&xy\ot r_x\ot 1  +x\ot r_y\ot xy,\\ t_1(yxy\ot x\ot 1)&=&
1\ot r_y\ot y+y\ot r_y \ot 1,\\

   t_1(xyxy\ot x\ot 1)&=&1\ot r_x\ot yxy-x\ot
r_x\ot x +yxy\ot r_x\ot 1 +yx\ot r_y\ot
xy\\

 t_1(1\ot y\ot 1)&=&0, \\
t_1(x\ot y\ot 1)&=&0,\\
t_1(y\ot y\ot 1)&=&1\ot r_y\ot 1,\\
 t_1(xy\ot y\ot 1)&=&1\ot r_x\ot yx +x\ot r_y\ot 1  +yx\ot r_y\ot x, \\
 t_1(yx\ot y\ot1)&=&0,\\
 t_1(xyx\ot y\ot 1)&=&0,\\
 t_1(yxy\ot y\ot 1)&=&y\ot r_x \ot yx+yx\ot r_y\ot 1,\\
  t_1(xyxy\ot y\ot 1)&=& xy\ot r_x\ot yx+xyx\ot r_y\ot 1.

\end{array}$$
Notice that  $t_1(b_1\ot b_2\ot 1)=0$ for $b_1, b_2\in
\mathcal{B}$ with $b_1b_2\in \mathcal{B}$. This observation will
simplify very much some computations.

The map $t_2: P_2\to P_3$ is given by
$$\begin{array}{lcl}
t_2(1\ot r_x\ot 1)&=&0, \\

  t_2(x\ot r_x\ot 1)&=&1 \ot1,\\

t_2(y\ot r_x\ot 1)&=&0,\\

t_2(xy\ot r_x\ot 1)&=&0, \\

 t_2(yx\ot r_x\ot 1)&=&y \ot 1,\\

 t_2(xyx\ot
r_x\ot 1)&=&xy\ot   1  +x \ot y,\\


t_2(yxy\ot r_x\ot 1)&=& 1\ot x,\\


t_2(xyxy\ot r_x\ot 1)&=& 1 \ot yxy  +yxy\ot 1+ y \ot  xy +yx\ot
y\\


 t_2(1\ot r_y\ot 1)&=&0, \\

t_2(x\ot r_y\ot 1)&=&0,\\

t_2(y\ot r_y\ot 1)&=&0,\\

t_2(xy\ot r_y\ot 1)&=&x \ot 1, \\


 t_2(yx\ot r_y\ot1)&=&0,\\

 t_2(xyx\ot r_y\ot 1)&=&0,\\

 t_2(yxy\ot r_y\ot 1)&=&y\ot x+yx\ot 1,\\

 t_2(xyxy\ot r_y\ot 1)&=& x\ot yx  +xy\ot  x+xyx\ot 1.


\end{array}$$

We define $\tau: P_3=A\ot A\to A$ as follows: $\tau(xyxy\ot 1)=1$
and $\tau(b\ot 1)=0$ for $b\in \mathcal{B}-\{xyxy\}$.  We impose
$t_3=t_{-1}\circ \tau: P_3\to P_4$ and define $t_{4n+i}=t_i$ for
$n\geq 0$ and $i\in \{0, 1, 2, 3\}$.

\begin{Prop} The above defined maps $\{t_i\}_{i\geq -1}$ form a weak self-homotopy
over $P_*$.

\end{Prop}

\begin{Proof}

Since the resolution is periodic of period $4$, it suffices to
prove that $$\left\{\begin{array}{lcl} d_0t_{-1}&=&Id, \\
 d_{p+1} t_p+t_{p-1} d_p&=&Id,\   \mathrm{for}\  0\leq p\leq 2,\\
  t_2d_3 +\rho
\tau&=&Id, \\   \tau\rho&=&Id.\end{array}\right.$$

The first two maps $t_{-1}$ and $t_0$ are given  at the end of
Section~\ref{Method}.

The map $t_1$  can be computed using the formula given at the end
of Section~\ref{Method}.  For instance, for $t_1(xyxy\ot x\ot 1)$,
one can write
$$xyxyx=x r_x x+ yxy r_x +1 r_x yxy +yx r_y xy + r_x^2 \in TV.$$
As $r_x^2\in I^2$, we have
$$t_1(xyxy\ot x\ot 1)=x\ot  r_x \ot x+ yxy\ot  r_x \ot 1+1\ot  r_x\ot  yxy +yx\ot  r_y \ot xy .$$
Another expression is
$$xyxyx=r_y yx +y r_y x+yxy r_x +yx r_y xy + yxx yxxy\in TV,$$
Notice that $yxx\in I$ and $yxxyxxy\in I^2$, which give
$$t_1(xyxy\ot x\ot 1)=1\ot r_y\ot  yx +y \ot r_y \ot x+yxy \ot r_x\ot 1 +yx \ot r_y\ot  xy. $$

\bigskip

The maps $t_2$ and $\tau$ are computed by direct inspection. The
details are tedious and long, but not difficult.

\end{Proof}

\section{Comparison morphisms for $kQ_8$}\label{Phi_Psi}

For an algebra $A$, denote by $\ol{A}=A/(k\cdot 1)$. The
normalized bar resolution is a quotient complex of the usual bar
resolution whose $p$-th term is   $B_p(A)=A\otimes \ol{A}^{\otimes
p}\otimes A$ and whose  differential is induced from that of the
usual bar resolution. It is easy to see that this complex is
well-defined.

Using the method from Section~\ref{Method}, one can compute
comparison morphsims between the minimal resolution $P_*$ and the
normalized bar resolution $\mathrm{Bar}_*(A)$, denoted by $\Phi_*: P_*\to
\mathrm{Bar}_*(A)$ and $\Psi_*: \mathrm{Bar}_*(A)\to P_*$.

The chain map $\Phi_*: P_*\to B_*:=\mathrm{Bar}_*(A)$  can be computed by
applying Lemma \ref{Phi}. Let us give the formulas for $\Phi_i $ with  $i\leq 5$.

 \begin{itemize}

\item  $\Phi_0=Id: P_0=A\ot A\to B_0=A\ot A$;

\item $\Phi_1: P_1=A\ot kQ_1\ot  A\to B_1=A\ot \ol{A}\ot  A$ is
induced by the inclusion $kQ_1\hookrightarrow  \ol{A}$;

\item $\Phi_2: P_2=A\ot kQ_1^*\ot  A\to B_2=A\ot \ol{A}^{\ot 2}\ot
A$ is given  by
$$\Phi_2(1\ot r_x\ot 1)=1\ot x\ot x\ot 1+1\ot y \ot x\ot y+1\ot
yx\ot y\ot 1$$ and
$$\Phi_2(1\ot r_y\ot 1)=1\ot y\ot y\ot 1+1\ot x \ot y\ot x+1\ot
xy\ot x\ot 1;$$

\item $\Phi_3: P_3=A\ot    A\to B_3=A\ot \ol{A}^{\ot 3}\ot A$ is
given  by
$$\begin{aligned}
\Phi_3(1\ot  1)=&1\ot x\ot x\ot x\ot 1+1\ot x\ot y \ot x\ot y+1\ot x\ot yx\ot y\ot 1\\
+&1\ot y\ot y\ot y\ot 1+1\ot y\ot x \ot y\ot x+1\ot y\ot xy\ot x\ot 1;\end{aligned}$$\label{Phi3}

\item $\Phi_4: P_4=A\ot    A\to B_4=A\ot \ol{A}^{\ot 4}\ot A$ is
given  by
$$\begin{aligned}
\Phi_4(1\ot  1)=\sum\limits_{b\in\mathcal{B}\setminus\{1\}}1\otimes b\Phi_3(1\otimes 1)b^*;\end{aligned}$$

\item $\Phi_5: P_5=A\ot  kQ_1\ot   A\to B_5=A\ot \ol{A}^{\ot 5}\ot
A$ is given  by
$$\begin{aligned}\Phi_5(1\ot x\ot  1)=1\otimes x\Phi_4(1\otimes 1)\end{aligned}$$
 and
 $$\begin{aligned}\Phi_5(1\ot y\ot  1)=1\otimes y\Phi_4(1\otimes 1).\end{aligned}$$
\end{itemize}

The chain map $\Psi_*: \mathrm{Bar}_*(A)\to P_*$ can be computed by
applying the method of Section~\ref{Method} to $t_*$. But the dimension of $\mathrm{Bar}_n(A)$ grows very fast. We have to specify the value of $\Psi_n$ on $7^n$ elements to fully describe it. So we give the full description only for $\Psi_0$ and $\Psi_1$.

\begin{itemize}
\item  $\Psi_0=Id: B_0=A\ot A\to P_0=A\ot A$;

\item $\Psi_1: B_1=A\ot \ol{A}\ot  A \to P_1=A\ot kQ_1\ot  A$ is
given  by  $\Psi_1(1\ot b\ot 1)={\sf C}(b)$ for $b\in
\mathcal{B}-\{1\}$.
\end{itemize}

\section{BV-structure on $HH^*(kQ_8)$}\label{BV}

Generalov proved the following result in \cite{Generalov}.

\begin{Thm}\label{Gen}\cite[Theorem 1.1, case  1b)]{Generalov}Let $k$ be an algebrically closed field of characteristic two. Let
$Q_8$ be the quaternion group of order $8$.  We have
$HH^*(kQ_8)\simeq k[\mathcal{X}]/I$ where
 \begin{itemize}\item $\mathcal{X}=\{p_1, p_2, p_2',
p_3, u_1, u_1', v_1, v_2, v_2', z\}$ with
$$\left\{\begin{array}{c} |p_1|=|p_2|=|p_2'|=|p_3|=0, |u_1|=|u_1'|=1,\\
 |v_1|=|v_2|=|v_2'|=2, |z|=4;\end{array}\right.$$

\item the ideal $I$ is generated by the following relations

\begin{itemize}
\item[] of degree $0$
$$\left\{\begin{array}{c} p_1^2, p_2^2, (p_1')^2, p_1p_2, p_1p_2',
p_2p_2',\\  p_3^2, p_1p_3, p_2p_3, p_2'p_3;\end{array}\right.$$

\item[]  of degree $1$
$$ \begin{array}{c} p_2u_1-p_2'u_1', p_2'u_1-p_1u_1', p_1u_1-p_2u_1';\end{array} $$

\item[] of degree $2$
$$\left\{\begin{array}{c}  p_1v_1, p_2v_2, p_2'v_2', p_3v_1, p_3v_2, p_3v_2', u_1u_1',\\
p_2v_1-p_1v_2', p_2v_1-p_2'v_2,
p_2v_1-p_3u_1^2,\\
p_2'v_1-p_1v_2, p_2'v_1-p_2v_2',
p_2'v_1-p_3(u_1')^2;\end{array}\right.$$

\item[] of degree $3$
$$ \begin{array}{c}  u_1'v_2-u_1v_2', u_1'v_1-u_1v_2, u_1v_1-u_1'v_2', u_1^3-(u_1')^3;\end{array} $$

\item[]of degree $4$
$$ \begin{array}{c}  v_1^2, v_2^2, (v_2')^2, v_1v_2, v_1v_2', v_2v_2'.\end{array}$$
\end{itemize}
\end{itemize}
\end{Thm}

\begin{Rem}\label{denot_coc_cob} Let $P$ be one of the members of the minimal resolution $P_*$. We use the following notion for the elements of $\mathrm{Hom}_{A^e}(P, A)$. If $P=A\ot A$ and $a\in A$, then we denote by $a$ the map which sends $1\otimes 1$ to  $a$.
If $P=A\ot Q_1\ot A$ ($P=A\ot Q_1^*\ot A$), $a,b\in A$, then we denote by $(a,b)$ the the map which sends $1\otimes x\otimes 1$ and $1\otimes y\otimes 1$ ($1\otimes r_x\otimes 1$ and $1\otimes r_y\otimes 1$) to $a$ and $b$ respectively. Moreover, we use the same notation for the corresponding cohomology classes. It follows from the work \cite{Generalov} that $p_1=xy+yx$, $p_2=xyx$, $p_2'=yxy$, $p_3=xyxy$,
$u_1=(1+xy, x)$, $u_1'=(y, 1+yx)$, $v_1=(y,x)$, $v_2=(x,0)$, $v_2'=(0,y)$ and $z=1$ in this notation. By \cite[Remarks 3.0.3, 3.1.18]{Generalov} we have that $$(xy+yx,0),(0,xy+yx),(xyx, yxy)\in B^1(A)$$
 and $$(xy+yx,yxy),(xyx,xy+yx),(xyxy,yxy),(xyx,xyxy)\in B^2(A).$$
\end{Rem}

We want to compute the Lie bracket and BV structure on $HH^*(kQ_8)$. By definition \ref{Def-BV-algebra} and the Poisson rule,
$$[a\smile b,\,  c]=[a,\, c]\smile b+(-1)^{|a|(|c|-1)}a\smile [b,\, c],$$
we have an equality (in characteristic 2)
\begin{equation}\label{BV-eq}
\Delta(abc)=\Delta(ab)c+\Delta(ac)b+\Delta(bc)a+\Delta(a)bc+\Delta(b)ac+\Delta(c)ab.
\end{equation}
So we need to compute $\Delta(x)$ only for $x\in \mathcal{X}$ and $x=a\cup b$ where $a, b\in \mathcal{X}$.
Suppose that $a\in HH^n(kQ_8)$ is given by a cocycle $f: P_n\to A$, then we compute $\Delta (a)$ using the following formula $$\Delta(a)=\Delta(f\circ  \Psi_n)\circ \Phi_{n-1}.$$
It is clear that $\Delta(a)=0$ for $a\in\{p_1,p_2,p_2',p_3\}$ because $\Delta$ is a map of degree $-1$.

For $b,c\in\mathcal{B}$ we have
$$\langle b, c\rangle=\begin{cases}1,&\mbox{if $c=b^*$,}\\0&\mbox{overwise.}\end{cases}$$
Then it follows from Theorem \ref{Tradler} that
$$
\Delta(\alpha)(a_1\otimes\dots\otimes a_{n-1})=\sum_{b\in\mathcal{B}\setminus\{1\}}\left\langle\sum_{i=1}^n (-1)^{i(n-1)}\alpha(a_i\otimes
\dots\otimes a_{n-1}\otimes b\otimes a_1\otimes \dots\otimes a_{i-1}), \; 1\right\rangle b^*
$$
for $\alpha\in C^n(A)$, $a_1, \dots, a_{n-1}\in A$.

\begin{Lem}\label{BV-structure1}
$$\begin{array}{ll}
\Delta(u_1)=\Delta(u_1')=0,&\Delta(p_1u_1)=\Delta(p_3u_1)=\Delta(p_2u_1')=p_2',\\
\Delta(p_2u_1)=\Delta(p_2'u_1')=p_1,&\Delta(p_2'u_1)=\Delta(p_1u_1')=\Delta(p_3u_1')=p_2.
\end{array}$$
\end{Lem}
\begin{Proof}  We have
$$
\begin{aligned}
p_1u_1&=p_2u_1'=(xyxy, xyx),\, p_2u_1=p_2'u_1'=(xyx, 0),\,p_2'u_1=p_1u_1'=(yxy, xyxy),\\
p_3u_1&=(xyxy, 0),\,p_3u_1=(0, xyxy)
\end{aligned}
$$
in $HH^1(A)$ (see Remark \ref{denot_coc_cob}).

For $a\in \mathrm{HH}^1(A)$ we have
$$
\Delta(a)(1\ot 1)=\Delta(a\circ \Psi_1)\Phi_0(1\ot 1)=\sum_{b\in\mathcal{B}\setminus\{1\}}\left\langle a\big({\sf C}(b)\big), \; 1\right\rangle b^*.
$$
It is easy to check that
$$\left\langle a\big({\sf C}(b)\big), \; 1\right\rangle=
\begin{cases}
0,&\mbox{if $a\in\{u_1,u_1'\}$, $b\in \mathcal{B}$, or $a\in\{p_1u_1,p_3u_1\}$, $b\in\mathcal{B}\setminus\{x\}$,}\\
&\mbox{or $a=p_2u_1$, $b\in\mathcal{B}\setminus\{xy,yx\}$, or $a\in\{p_2'u_1,p_3u_1'\}$, $b\in\mathcal{B}\setminus\{y\}$,}\\
1,&\mbox{if $a\in\{p_1u_1,p_3u_1\}$, $b=x$, or $a=p_2u_1$, $b\in\{xy,yx\}$,}\\
&\mbox{or $a\in\{p_2'u_1,p_3u_1'\}$, $b=y$.}\\
\end{cases}$$
Lemma follows from this formula.
\end{Proof}

\begin{Lem}\label{BV-structure2}
$$\Delta(ab)=0.$$
for $a\in \{v_1, v_2, v_2'\}, b\in \{1,p_1,p_2,p_2',p_3\}$.
\end{Lem}
\begin{Proof} For $a\in \mathrm{HH}^2(A)$ we have
$$
\begin{aligned}
\Delta(a)(1\ot x\ot 1)&=\Delta(a\circ \Psi_2)\Phi_1(1\ot x\ot 1)=\sum_{b\in\mathcal{B}\setminus\{1\}}\left\langle(a\circ \Psi_2)(b\ot x+x\ot b), \; 1\right\rangle b^*,\\
\Delta(a)(1\ot y\ot 1)&=\Delta(a\circ \Psi_2)\Phi_1(1\ot y\ot 1)=\sum_{b\in\mathcal{B}\setminus\{1\}}\left\langle(a\circ \Psi_2)(b\ot y+y\ot b), \; 1\right\rangle b^*.
\end{aligned}
$$

Direct calculations show that
$$\begin{aligned}
&\Psi_2(b\ot x+x\ot b)=t_1(b\ot x\ot 1+x{\sf C}(b))\\
=&\begin{cases}
0,&\mbox{if $b\in\{x,y\}$,}\\
+1\ot r_x\ot y+y\ot r_x\ot yx+yx\ot r_y\ot 1,&\mbox{if $b=xy$,}\\
y\ot r_x\ot 1+xy\ot r_x\ot y+1\ot r_y\ot xy,&\mbox{if $b=yx$,}\\
xy\ot r_x\ot 1+x\ot r_y\ot xy+1\ot r_x\ot yx+yx\ot r_y\ot x,&\mbox{if $b=xyx$,}\\
1\ot r_y\ot y+y\ot r_y\ot 1,&\mbox{if $b=yxy$,}\\
x\ot r_x\ot x+yxy\ot r_x\ot 1,&\mbox{if $b=xyxy$;}\\
\end{cases}
\end{aligned}
$$
$$\begin{aligned}
&\Psi_2(b\ot y+y\ot b)=t_1(b\ot y\ot 1+y{\sf C}(b))\\
=&\begin{cases}
0,&\mbox{if $b\in\{x,y\}$,}\\
1\ot r_x\ot yx+x\ot r_y\ot 1+yx\ot r_y\ot x,&\mbox{if $b=xy$,}\\
+1\ot r_y\ot x+xy\ot r_x\ot 1+x\ot r_y\ot xy,&\mbox{if $b=yx$,}\\
1\ot r_y\ot y+y\ot r_y\ot 1,&\mbox{if $b=xyx$,}\\
y\ot r_x\ot yx+yx\ot r_y\ot 1+1\ot r_y\ot xy+xy\ot r_x\ot y,&\mbox{if $b=yxy$,}\\
1\ot r_y\ot xyx+y\ot r_y\ot y,&\mbox{if $b=xyxy$.}\\
\end{cases}
\end{aligned}
$$
It follows from Remark \ref{denot_coc_cob} and the formulas above that $\Delta(v_1)=\Delta(v_2)=\Delta(v_2')=0$ in $HH^1(A)$.

The remaining formulas of lemma can be deduced in the same way. But there is an easier way. By Theorem \ref{Gen} it is enough to prove that $\Delta(p_3u_1^2)=\Delta(p_3(u_1')^2)=0$. And this equalities can be easily deduced from Lemma \ref{BV-structure1} and the formula \eqref{BV-eq}.
\end{Proof}

\begin{Lem}\label{BV-structure3}
$$\Delta(u_1v_1)=\Delta(u_1^{\prime}v_2^{\prime})=(u_1^{\prime})^{2}+v_2,\Delta(u_1^{\prime}v_1)=\Delta(u_1v_2)=u_1^{2}+v_2^{\prime}, \Delta(u_1^{\prime}v_2)=\Delta(u_1v_2')=v_1$$
in $HH^2(A).$
\end{Lem}
\begin{Proof} For $a\in \mathrm{HH}^3(A)$ we have
$$
\begin{aligned}
\Delta(a)(1\ot r_x\ot 1)&=\Delta(a\circ \Psi_3)\Phi_2(1\ot r_x\ot 1)\\
&=\sum_{b\in\mathcal{B}\setminus\{1\}}\left\langle(a\circ \Psi_3)
(b\ot x\ot x+x\ot b\ot x+x\ot x\ot b), \; 1\right\rangle b^*\\
&+\sum_{b\in\mathcal{B}\setminus\{1\}}\left\langle(a\circ \Psi_3)
(b\ot y\ot x+x\ot b\ot y+y\ot x\ot b), \; 1\right\rangle b^*y\\
&+\sum_{b\in\mathcal{B}\setminus\{1\}}\left\langle(a\circ \Psi_3)
(b\ot yx\ot y+y\ot b\ot yx+yx\ot y\ot b), \; 1\right\rangle b^*,\\
\Delta(a)(1\ot r_y\ot 1)&=\Delta(a\circ \Psi_3)\Phi_2(1\ot r_y\ot 1)\\
&=\sum_{b\in\mathcal{B}\setminus\{1\}}\left\langle(a\circ \Psi_3)
(b\ot y\ot y+y\ot b\ot y+y\ot y\ot b), \; 1\right\rangle b^*\\
&+\sum_{b\in\mathcal{B}\setminus\{1\}}\left\langle(a\circ \Psi_3)
(b\ot x\ot y+y\ot b\ot x+x\ot y\ot b), \; 1\right\rangle b^*x\\
&+\sum_{b\in\mathcal{B}\setminus\{1\}}\left\langle(a\circ \Psi_3)
(b\ot xy\ot x+x\ot b\ot xy+xy\ot x\ot b), \; 1\right\rangle b^*.
\end{aligned}
$$

Direct calculations (see also the proof of Lemma \ref{BV-structure2}) show that
$$\begin{aligned}
&\Psi_3(b\ot x\ot x+x\ot b\ot x+x\ot x\ot b)=t_2(b\ot r_x\ot 1+xt_1(b\ot x\ot 1+x{\sf C}(b)))\\
=&\begin{cases}
1\ot 1,&\mbox{if $b=x$,}\\
0,&\mbox{if $b=y$,}\\
1\ot y,&\mbox{if $b=xy$,}\\
y\ot 1,&\mbox{if $b=yx$,}\\
xy\ot 1+x\ot y+yx\ot xy+1\ot yx,&\mbox{if $b=xyx$,}\\
1\ot x+x\ot 1,&\mbox{if $b=yxy$,}\\
1\ot yxy,&\mbox{if $b=xyxy$;}
\end{cases}\\
&\Psi_3(b\ot y\ot x+x\ot b\ot y+y\ot x\ot b)=t_2(xt_1(b\ot y\ot 1)+yt_1(x{\sf C}(b)))\\
=&\begin{cases}
0,&\mbox{if $b\not=xy$,}\\
1\ot yx+y\ot x+yx\ot 1+xy\ot yx,&\mbox{if $b=xy$;}\\
\end{cases}\\
 \end{aligned}
$$
 $$
\begin{aligned}
&\Psi_3(b\ot yx\ot y+y\ot b\ot yx+yx\ot y\ot b)=t_2(yt_1(b\ot y\ot x+by\ot x\ot 1)+yxt_1(y{\sf C}(b)))\\
=&\begin{cases}
0,&\mbox{if $b\in\{x,xy,yx,yxy\}$,}\\
1\ot x,&\mbox{if $b=y$,}\\
yx\ot 1,&\mbox{if $b=xyx$,}\\
xy\ot yxy+xyx\ot x,&\mbox{if $b=xyxy$;}
\end{cases}\\
&\Psi_3(b\ot y\ot y+y\ot b\ot y+y\ot y\ot b)=t_2(b\ot r_y\ot 1+yt_1(b\ot y\ot 1+y{\sf C}(b)))\\
=&\begin{cases}
0,&\mbox{if $b\in\{x,y,xyx\}$,}\\
x\ot 1,&\mbox{if $b=xy$,}\\
1\ot x,&\mbox{if $b=yx$,}\\
y\ot x+yx\ot 1+xy\ot yx+1\ot xy,&\mbox{if $b=yxy$,}\\
x\ot yx+xy\ot x+xyx\ot 1,&\mbox{if $b=xyxy$;}\\
\end{cases}\\
&\Psi_3(b\ot x\ot y+y\ot b\ot x+x\ot y\ot b)=t_2(yt_1(b\ot x\ot 1)+xt_1(y{\sf C}(b)))\\
=&\begin{cases}
0,&\mbox{if $b\in\{x,y,xy,yxy\}$,}\\
xy\ot 1+x\ot y+1\ot xy+yx\ot xy,&\mbox{if $b=yx$,}\\
1\ot x+x\ot 1,&\mbox{if $b=xyx$,}\\
y\ot x,&\mbox{if $b=xyxy$;}
\end{cases}\\
&\Psi_3(b\ot xy\ot x+x\ot b\ot xy+xy\ot x\ot b)=t_2(xt_1(b\ot x\ot y+bx\ot y\ot 1)+xyt_1(x{\sf C}(b)))\\
=&\begin{cases}
1\ot y,&\mbox{if $b=x$,}\\
0,&\mbox{if $b\in\{y,xy,yx,xyx\}$,}\\
x\ot y,&\mbox{if $b=yxy$,}\\
yxy\ot y+yx\ot xyx,&\mbox{if $b=xyxy$.}
\end{cases}
\end{aligned}
$$

By Theorem \ref{Gen} it is enough to calculate $\Delta$ on $u_1'v_2'$, $u_1v_2$ and $u_1'v_2$. By \cite[Lemmas 4.1.2, 4.1.8]{Generalov} and Remark \ref{denot_coc_cob} we have
$$
\begin{aligned}
u_1'v_2'&=u_1'T^1(v_2')=(y,1+yx)\begin{pmatrix} 0\\y\ot 1\end{pmatrix}=y,\:u_1v_2=u_1T^1(v_2)=(1+xy,x)\begin{pmatrix} x\ot 1\\0\end{pmatrix}=x,\\
u_1'v_2&=u_1T^1(v_2)=(y,1+yx)\begin{pmatrix} x\ot 1\\0\end{pmatrix}=xy,\\
u_1^2&=u_1\T^1(u_1)=(1+xy,x)\begin{pmatrix} 1\ot(1+xy+yx) & (y+yxy)\ot 1\\
  1\ot y + x\ot x + yxy\ot 1& x\ot 1 + 1\ot x + x\ot yx
 \end{pmatrix}=(1,y),\\
 (u_1')^2&=u_1'\T^1(u_1')=(y,1+yx)\begin{pmatrix} y\ot 1+1\ot y+y\ot xy & 1\ot x+y\ot y+xyx\ot 1\\
 (x+xyx)\ot 1& 1\ot(1+yx+xy)
 \end{pmatrix}=(x,1).\\
\end{aligned}
$$
From the formulas above we obtain
$$
\Delta(u_1'v_2')=(0,1)=(u_1')^2+v_2,\,\Delta(u_1v_2)=(1,0)=u_1^2+v_2',\,\Delta(u_1'v_2)=(y,x)=v_1.
$$
\end{Proof}

\begin{Lem}\label{BV-structure4}
$$\Delta(az)=0.$$
for $a\in\{a,p_1,p_2,p_2',p_3\}$.
\end{Lem}

\begin{Proof} It follows from the formula for $\Phi_3$ that we have to calculate $\Psi_4$ on four kinds of elements:\\
1) $b\ot a_1\ot a_2\ot a_3$;\\
2) $a_3\ot b\ot a_1\ot a_2$;\\
3) $a_2\ot a_3\ot b\ot a_1$;\\
4) $a_1\ot a_2\ot a_3\ot b$.\\
In all points $b\in \mathcal{B}$,
$$(a_1,a_2,a_3)\in \mathcal{A}=\{(x,x,x),(x,y,x),(x,yx,y),(y,y,y),(y,x,y),(y,xy,x)\}.$$

1) Note that
$$\Psi_3(a_1\ot a_2\ot a_3)=t_2(a_1t_1(a_2\otimes a_3\otimes 1))=\begin{cases}1\ot 1,&\mbox{if $a_1=a_2=a_3=x$,}\\0&\mbox{if $(a_1,a_2,a_3)\in \mathcal{A}\setminus (x,x,x)$.}\end{cases}$$
So if $(a_1,a_2,a_3)\in \mathcal{A}$, then
$$\Psi_4(b\ot a_1\ot a_2\ot a_3)=t_3(b\Psi_3(a_1\ot a_2\ot a_3))=\begin{cases}1\otimes 1,&\mbox{if $b=xyxy$, $a_1=a_2=a_3=x$,}\\0&\mbox{overwise.}\end{cases}$$

2) If $(a_1,a_2,a_3)\in\mathcal{A}\setminus\{(x,x,x),(y,y,y)\}$, then $t_1(a_1{\sf C}(a_2))=0$ and so $\Psi_4(a_3\ot b\ot a_1\ot a_2)=0$. For the remaining cases we have
$$
\begin{aligned}
\Psi_4(x\ot b\ot x\ot x)&=t_3(xt_2(b\otimes r_x\otimes 1))=
\begin{cases}1\otimes 1,&\mbox{if $b=xyxy$,}\\0&\mbox{overwise;}\end{cases}\\
\Psi_4(y\ot b\ot y\ot y)&=t_3(yt_2(b\otimes r\otimes 1))=
\begin{cases}1\otimes 1,&\mbox{if $b=xyxy$,}\\0&\mbox{overwise.}\end{cases}
\end{aligned}
$$

Let $b\in\mathcal{B}$, $r\in\{r_x,r_y\}$. Note that
$t_3(xt_2(b\ot r\ot 1))$ can be nonzero only for $(b,r)=(xyxy, r_x)$. Analogously $t_3(yt_2(b\ot r\ot 1))$ can be nonzero only for $(b,r)=(xyxy, r_y)$.
Also note that for $b\in\mathcal B$, $a\in\{x,y\}$ the element $t_1(b\otimes a\otimes 1)$ is a sum of elements of the form $u\otimes r\otimes v$, where $u,v\in \mathcal B$, $r\in\{r_x,r_y\}$ and $(u,r)\not\in\{(xyx,r_x),(yx,r_x),(yxy,r_y),(xy,r_y)\}$. So we have the equalities
$$
t_3(xt_2(yt_1(A))=t_3(xt_2(yxt_1(A))=t_3(yt_2(xt_1(A))=t_3(yt_2(xyt_1(A))=0
$$
for any $A\in P_1$. In the same way the equalities
$$
t_3(yxt_2(yt_1(A))=t_3(xyt_2(xt_1(A))=0
$$
can be proved. Then $\Psi_4$ can be nonzero in points 3) and 4) only for $a_1=a_2=a_3=x$ and $a_1=a_2=a_3=y$.
The same arguments show that
$$t_3(xt_2(xt_1(b\otimes a\otimes 1))=0\,\:((b,a)\in (\mathcal{B}\times\{x,y\})\setminus\{xyxy,x\})$$
and
$$t_3(yt_2(yt_1(b\otimes a\otimes 1))=0\,\:((b,a)\in (\mathcal{B}\times\{x,y\})\setminus\{xyxy,y\}).$$
So we obtain equalities
$$
\begin{aligned}
\Psi_4(x\ot x\ot b\ot x)&=\begin{cases}1\otimes 1,&\mbox{if $b=xyxy$,}\\0&\mbox{overwise;}\end{cases}\,
\Psi_4(x\ot x\ot x\ot b)=0;\\
\Psi_4(y\ot y\ot b\ot y)&=\Psi_4(y\ot y\ot y\ot b)=\begin{cases}1\otimes 1,&\mbox{if $b=xyxy$,}\\0&\mbox{overwise.}\end{cases}
\end{aligned}
$$
We set
$$
\begin{aligned}
S(a_1,a_2,a_3,b)&:=\Psi_4(b\ot a_1\ot a_2\ot a_3+a_3\ot b\ot a_1\ot a_2+a_2\ot a_3\ot b\ot a_1+a_1\ot a_2\ot a_3\ot b)\\
&=\begin{cases}1\otimes 1,&\mbox{if $b=xyxy$, $(a_1,a_2,a_3)\in\{(x,x,x),(y,y,y)\}$,}\\0&\mbox{if $b\in\mathcal{B}\setminus\{xyxy\}$ or $(a_1,a_2,a_3)\in\mathcal{A}\setminus\{(x,x,x),(y,y,y)\}$.}\end{cases}
\end{aligned}
$$
Then for $a\in HH^4(A)$ we have
$$
\begin{aligned}
\Delta(a)(1\ot 1)&=\Delta(a\circ \Psi_4)\Phi_3(1\ot 1)\\
&=\sum_{b\in\mathcal{B}\setminus\{1\},(a_1,a_2,a_3)\in\mathcal{A}\setminus\{(x,y,x),(y,x,y)\}}\left\langle a(S(a_1,a_2,a_3,b)), \; 1\right\rangle b^*\\
&+\sum_{b\in\mathcal{B}\setminus\{1\}}\left\langle a(S(x,y,x,b)), \; 1\right\rangle b^*y+\sum_{b\in\mathcal{B}\setminus\{1\}}\left\langle a(S(y,x,y,b)), \; 1\right\rangle b^*x\\
&=\left\langle a(1\ot 1+1\ot 1), \; 1\right\rangle=0.
\end{aligned}
$$
\end{Proof}

If we know the values of $\Delta(a)$ and $\Delta(b)$, then it is enough to calculate $[a,b]$ to find $\Delta(ab)$. Sometimes it is easier than calculate $\Delta(ab)$ directly. Suppose that $a$ and $b$ are given by cocycles $f: P_n\to A$ and $g: P_m\to A$, then we compute $[a, b]$ using the following formula
$$[a, b]=   [f\circ \Psi_n, g\circ \Psi_m]\circ \Phi_{n+m-1}. $$

\begin{Lem}\label{BV-structure5}
$$\Delta(u_1z)=\Delta(u_1'z)=0.$$
\end{Lem}
\begin{Proof}
It is enough to prove that $[u_1,z]=[u_1',z]=0$. For $a\in\{u_1,u_1'\}$ we have
$$
[a,z](1\ot 1)=\big((a\circ \Psi_1)\circ (z\circ \Psi_4)\big)\Phi_4(1\ot 1)+\big((z\circ \Psi_4)\circ(a\circ \Psi_1)\big) \Phi_4(1\ot 1).
$$
Let prove that $\Psi_4\Phi_4=Id$. Direct calculations show that $\Psi_3\Phi_3=Id$ (see the proof of Lemma \ref{BV-structure4}). Then
$$
(\Psi_4\Phi_4)(1\ot 1)=\sum\limits_{b\in\mathcal{B}\setminus\{1\}}t_3(b\Psi_3\Phi_3(1\otimes 1))b^*=\sum\limits_{b\in\mathcal{B}\setminus\{1\}}t_3(b\otimes 1)b^*=1\ot 1.
$$
If $|a|=1$, we have
$$
\big((a\circ \Psi_1)\circ (z\circ \Psi_4)\big)\Phi_4(1\ot 1)=(a\circ \Psi_1)\big(z(\Psi_4\Phi_4)(1\ot 1)\big)=(a\circ \Psi_1)(1)=0.
$$
It remains to prove that
$$
\big((z\circ \Psi_4)\circ(a\circ \Psi_1)\big) \Phi_4(1\ot 1)=0
$$
for $a\in\{u_1,u_1'\}$. Let us introduce the notation
$$f=(z\circ \Psi_4)\circ(u_1\circ \Psi_1),\,f'=(z\circ \Psi_4)\circ(u_1'\circ \Psi_1).$$

In this proof we need to know the values of $u_1\circ \Psi_1$ and $u_1'\circ \Psi_1$ on elements of $\mathcal{B}$. Direct calculations show that
$$
(u_1\circ \Psi_1)(b)=\begin{cases}
1+xy,&\mbox{if $b=x$},\\
x,&\mbox{if $b=y$},\\
y+yxy,&\mbox{if $b=xy$},\\
y,&\mbox{if $b=yx$},\\
xy+yx,&\mbox{if $b=xyx$},\\
xyx,&\mbox{if $b=yxy$},\\
yxy,&\mbox{if $b=xyxy$};
\end{cases}
(u_1'\circ \Psi_1)(b)=\begin{cases}
y,&\mbox{if $b=x$},\\
1+yx,&\mbox{if $b=y$},\\
x,&\mbox{if $b=xy$},\\
x+xyx,&\mbox{if $b=yx$},\\
yxy,&\mbox{if $b=xyx$},\\
xy+yx,&\mbox{if $b=yxy$},\\
xyx,&\mbox{if $b=xyxy$}.
\end{cases}
$$

We want to calculate the values of $(z\circ \Psi_4)\circ(a\circ \Psi_1)$ on elements of the form $b\otimes a_1\ot a_2\ot a_3$, where $b\in\mathcal{B}\setminus\{1\}$ and $(a_1,a_2,a_3)\in\mathcal{A}$ (see the proof of Lemma \ref{BV-structure4} for notation). Let consider each element of $\mathcal{A}$ separately.

1) Let $a_1=a_2=a_3=x$. Let $a\in\{u_1,u_1'\}$. Then
$$(z\circ \Psi_4)\circ_1(a\circ \Psi_1)(b\otimes x\otimes x\otimes x)=zt_3((a\circ \Psi_1)(b)\ot 1)=0,$$
because $(a\circ \Psi_1)(b)$ is a sum of elements of $\mathcal{B}\setminus\{xyxy\}$. Further we have
$$
\begin{aligned}
\Psi_3((u_1\circ \Psi_1)(x)\otimes x\otimes x)&=\Psi_3(xy\otimes x\otimes x)=t_2(xy\ot r_x\ot 1)=0;\\
\Psi_2((u_1\circ \Psi_1)(x)\otimes x)&=t_1(xy\ot x\ot 1)=0;\\
\Psi_3(x\otimes x\ot (u_1\circ \Psi_1)(x))&=\Psi_3(x\otimes x\otimes xy)=t_2(xt_1(x\ot x\ot y+yxy\ot y\ot 1))=1\ot y;\\
\Psi_3(y\otimes x\otimes x)&=\Psi_3(x\otimes y\otimes x)=\Psi_3(x\otimes x\otimes y)=0.\\
\end{aligned}
$$
Consequently,
$$
f(b\otimes x\otimes x\otimes x)=\begin{cases}
y,&\mbox{if $b=xyxy$},\\
0&\mbox{overwise,}
\end{cases}\,
f'(b\otimes x\otimes x\otimes x)=0.
$$

2) $(a_1,a_2,a_3)=(x,y,x)$. We have
$$(z\circ \Psi_4)\circ_1(a\circ \Psi_1)(b\otimes x\otimes y\otimes x)=(z\circ \Psi_4)\circ_2(a\circ \Psi_1)(b\otimes x\otimes y\otimes x)=0$$
for $a\in\{u_1,u_1'\}$, because $\Psi_2(y\ot x)=0$. Further we have
$$
\begin{aligned}
\Psi_3(x\ot (u_1\circ \Psi_1)(y)\otimes x)&=\Psi_3(x\ot x\otimes x)=1\ot 1;\\
\Psi_2(y\ot (u_1\circ \Psi_1)(x))&=\Psi_2(y\ot xy)=t_1(y\ot x\ot y+yx\ot y\ot 1)=0;\\
\Psi_3(x\ot (u_1'\circ \Psi_1)(y)\otimes x)&=\Psi_3(x\ot yx\otimes x)=t_2(xt_1(yx\ot x))\\
&=t_2(xy\ot r_x\ot 1+x\ot r_y\ot xy)=0;\\
\Psi_3(x\ot y\otimes (u_1'\circ \Psi_1)(x))&=\Psi_3(x\ot y\otimes y)=t_2(x\ot r_y\ot 1)=0.
\end{aligned}
$$
Consequently,
$$
f(b\otimes x\otimes y\otimes x)=\begin{cases}
1,&\mbox{if $b=xyxy$},\\
0&\mbox{overwise,}
\end{cases}\,
f'(b\otimes x\otimes y\otimes x)=0.
$$

3) $(a_1,a_2,a_3)=(x,yx,y)$. We have
$$(z\circ \Psi_4)\circ_1(a\circ \Psi_1)(b\otimes x\otimes yx\otimes y)=(z\circ \Psi_4)\circ_2(a\circ \Psi_1)(b\otimes x\otimes yx\otimes y)=0$$
for $a\in\{u_1,u_1'\}$, because $\Psi_2(yx\ot y)=0$. Further we have
$$
\begin{aligned}
\Psi_3(x\ot (u_1\circ \Psi_1)(yx)\otimes y)&=\Psi_3(x\ot y\otimes y)=t_2(x\ot r_y\ot 1)=0;\\
\Psi_3(x\ot yx\otimes (u_1\circ \Psi_1)(y))&=\Psi_3(x\ot yx\otimes x)=t_2(xy\ot r_x\ot 1+x\ot r_y\ot xy)=0;\\
\Psi_2((u_1'\circ \Psi_1)(yx)\otimes y)&=\Psi_2((x+xyx)\otimes y)=t_1(x\ot y\ot 1+xyx\ot y\ot 1)=0;\\
\Psi_3(x\ot yx\otimes (u_1\circ \Psi_1)(y))&=\Psi_3(x\ot yx\otimes yx)=t_2(xt_1(yx\ot y\ot x+yxy\ot x\ot 1))\\
&=t_2(x\ot r_y\ot y+xy\ot r_y\ot 1)=x\ot 1.\\
\end{aligned}
$$
Consequently,
$$
f(b\otimes x\otimes yx\otimes y)=0,\,f'(b\otimes x\otimes yx\otimes y)=\begin{cases}
1,&\mbox{if $b=yxy$},\\
0&\mbox{overwise.}
\end{cases}
$$

4) Let $a_1=a_2=a_3=y$. Let $a\in\{u_1,u_1'\}$. Then
$$(z\circ \Psi_4)\circ_1(a\circ \Psi_1)(b\otimes y\otimes y\otimes y)=0,$$
because $\Psi_3(y\ot y\ot y)=0$. Further we have
$$
\begin{aligned}
\Psi_3(x\otimes y\otimes y)&=\Psi_3(y\otimes x\otimes y)=\Psi_3(y\otimes y\otimes x)=0,\\
\Psi_3((u_1'\circ \Psi_1)(y)\otimes y\otimes y)&=\Psi_3(yx\otimes y\otimes y)=t_2(yx\ot r_y\ot 1)=0;\\
\Psi_2((u_1'\circ \Psi_1)(y)\otimes y)&=t_1(yx\ot y\ot 1)=0;\\
\Psi_3(y\otimes y\ot (u_1'\circ \Psi_1)(y))&=\Psi_3(y\otimes y\otimes yx)=t_2(yt_1(y\ot y\ot x+xyx\ot x\ot 1))=1\ot x.\\
\end{aligned}
$$
Consequently,
$$
f(b\otimes y\otimes y\otimes y)=0,\,f'(b\otimes y\otimes y\otimes y)=\begin{cases}
x,&\mbox{if $b=xyxy$},\\
0&\mbox{overwise.}
\end{cases}.
$$

5) $(a_1,a_2,a_3)=(y,x,y)$. We have
$$(z\circ \Psi_4)\circ_1(a\circ \Psi_1)(b\otimes y\otimes x\otimes y)=(z\circ \Psi_4)\circ_2(a\circ \Psi_1)(b\otimes y\otimes x\otimes y)=0$$
for $a\in\{u_1,u_1'\}$, because $\Psi_2(x\ot y)=0$. Further we have
$$
\begin{aligned}
\Psi_3(y\ot (u_1\circ \Psi_1)(x)\otimes y)&=t_2(yt_1(xy\otimes y\ot 1))=t_2(y\ot r_x\ot yx+yx\ot r_y\ot 1)=0;\\
\Psi_3(y\ot x\ot (u_1\circ \Psi_1)(y))&=\Psi_3(y\ot x\ot x)=0;\\
\Psi_3(y\ot (u_1'\circ \Psi_1)(x)\otimes y)&=\Psi_3(y\ot y\otimes y)=0;\\
\Psi_2(x\otimes (u_1'\circ \Psi_1)(y))&=\Psi_2(x\otimes yx)=t_1(x\ot y\ot x+xy\ot x\ot 1)=0.
\end{aligned}
$$
Consequently,
$$
f(b\otimes y\otimes x\otimes y)=f'(b\otimes y\otimes x\otimes y)=0.
$$

6) $(a_1,a_2,a_3)=(y,xy,x)$. We have
$$(z\circ \Psi_4)\circ_1(a\circ \Psi_1)(b\otimes y\otimes xy\otimes x)=(z\circ \Psi_4)\circ_2(a\circ \Psi_1)(b\otimes y\otimes xy\otimes x)=0$$
for $a\in\{u_1,u_1'\}$, because $\Psi_2(xy\ot x)=0$. Further we have
$$
\begin{aligned}
\Psi_3(y\ot (u_1\circ \Psi_1)(xy)\otimes x)&=\Psi_3(y\ot (y+yxy)\otimes x)=t_2(y\ot r_y\ot y-xyx\ot r_y\ot 1)=0;\\
\Psi_2(xy\otimes (u_1\circ \Psi_1)(x))&=\Psi_2(xy\otimes xy)=t_1(xy\ot x\ot y+xyx\ot y\ot 1)=0;\\
\Psi_3(y\ot (u_1'\circ \Psi_1)(xy)\otimes x)&=\Psi_3(y\ot x\otimes x)=0;\\
\Psi_3(y\ot xy\otimes (u_1'\circ \Psi_1)(x))&=\Psi_3(y\ot xy\otimes y)=t_2(y\ot r_x\ot yx+yx\ot r_y\ot 1)=0.\\
\end{aligned}
$$
Consequently,
$$
f(b\otimes y\otimes xy\otimes x)=f'(b\otimes y\otimes xy\otimes x)=0.
$$

Thus we obtain
$$
\begin{aligned}
f\Phi_4(1\ot 1)&=\sum\limits_{b\in\mathcal{B}\setminus\{1\}}f(b\ot x\ot x\ot x+b\ot x\ot yx\ot y+b\ot y\ot y\ot y+b\ot y\ot xy\ot x)b^*\\
&+\sum\limits_{b\in\mathcal{B}\setminus\{1\}}f(b\ot x\ot y\ot x)yb^*+\sum\limits_{b\in\mathcal{B}\setminus\{1\}}f(b\ot y\ot x\ot y)xb^*=y+y=0;\\
f'\Phi_4(1\ot 1)&=\sum\limits_{b\in\mathcal{B}\setminus\{1\}}f'(b\ot x\ot x\ot x+b\ot x\ot yx\ot y+b\ot y\ot y\ot y+b\ot y\ot xy\ot x)b^*\\
&+\sum\limits_{b\in\mathcal{B}\setminus\{1\}}f'(b\ot x\ot y\ot x)yb^*+\sum\limits_{b\in\mathcal{B}\setminus\{1\}}f'(b\ot y\ot x\ot y)xb^*=x+x=0;\\
\end{aligned}
$$

\end{Proof}

\begin{Lem}\label{BV-structure6}
$$\Delta(v_1z)=\Delta(v_2z)=\Delta(v_2'z)=0.$$
\end{Lem}
\begin{Proof} Firstly note that it is enough to prove that $[v_2,z]=0$. Indeed, by Jacoby identity and Lemmas \ref{BV-structure1}--\ref{BV-structure5} we have
$$
\begin{aligned}
\Delta(v_1z)&=[v_1,z]=[[u_1',v_2],z]=[v_2,[u_1',z]]+[u_1',[v_2,z]]=[u_1',[v_2,z]],\\
\Delta(v_2'z)&=[v_2',z]=[[u_1,v_2]+u_1^2,z]=[v_2,[u_1,z]]+[u_1,[v_2,z]]+2[u_1,z]=[u_1,[v_2,z]]
\end{aligned}
$$
and $\Delta(v_2z)=[v_2,z]$. For $a\in\{x,y\}$ we have
$$
[v_2,z](1\ot a\ot 1)=\big((v_2\circ \Psi_2)\circ (z\circ \Psi_4)\big)\Phi_5(1\ot a\ot 1)+\big((z\circ \Psi_4)\circ(v_2\circ \Psi_2)\big) \Phi_5(1\ot a\ot 1).
$$
Note that if $(a_1,a_2,a_3)\in \mathcal{A}\setminus\{(x,x,x),(y,y,y)\}$, then $\Psi_2(1\ot a_1\ot a_2\ot 1)=\Psi_2(1\ot a_2\ot a_3\ot 1)=0$ (see the proof of Lemma \ref{BV-structure4} for the notion of $\mathcal{A}$). It follows from this that
$$
((v_2\circ \Psi_2)\circ_i (z\circ \Psi_4)\big)(a\ot b\ot a_1\ot a_2\ot a_3)=0
$$
for $i\in\{1,2\}$, $a\in\{x,y\}$, $b\in\mathcal{B}$, $(a_1,a_2,a_3)\in \mathcal{A}\setminus\{(x,x,x),(y,y,y)\}$. Further we have
$$
\begin{aligned}
(z\circ \Psi_4)(b\ot y\ot y\ot y)&=zt_3(b\Psi_3(y\ot y\ot y))=0,\\
(z\circ \Psi_4)(a\ot b\ot y\ot y)&=zt_3(at_2(b\ot r_y\ot 1))=\begin{cases}
1,&\mbox{if $a=y$, $b=xyxy$,}\\
0,&\mbox{if $a=x$ or $b\in\mathcal{B}\setminus\{xyxy\}$.}\\
\end{cases}
\end{aligned}
$$
Because of $\Psi_2(1\ot 1\ot y\ot 1)=0$, we have
$$
\big((v_2\circ \Psi_2)\circ (z\circ \Psi_4)\big)(a\ot b\ot y\ot y\ot y)=0.
$$
Further we have
$$
\begin{aligned}
(z\circ \Psi_4)(b\ot x\ot x\ot x)&=zt_3(b\ot 1)=\begin{cases}
1,&\mbox{if $b=xyxy$,}\\
0,&\mbox{if $b\in\mathcal{B}\setminus\{xyxy\}$;}
\end{cases}\\
(z\circ \Psi_4)(a\ot b\ot x\ot x)&=zt_3(at_2(b\ot r_x\ot 1))=\begin{cases}
1,&\mbox{if $a=x$, $b=xyxy$,}\\
0,&\mbox{if $a=y$ or $b\in\mathcal{B}\setminus\{xyxy\}$.}\\
\end{cases}
\end{aligned}
$$
Because of $\Psi_2(1\ot a\ot 1\ot 1)=\Psi_2(1\ot 1\ot x\ot 1)=0$, we have
$$
\big((v_2\circ \Psi_2)\circ (z\circ \Psi_4)\big)(a\ot b\ot x\ot x\ot x)=0.
$$

It remains to prove that
$$
\big((z\circ \Psi_4)\circ(v_2\circ \Psi_2)\big) \Phi_5=0
$$
in $HH^5(A)$. If $(a_1,a_2,a_3)\in \mathcal{A}\setminus\{(x,x,x),(y,y,y)\}$, then
$$
\big((z\circ \Psi_4)\circ_i(v_2\circ \Psi_2)\big) (a\ot b\ot a_1\ot a_2\ot a_3)=0
$$
for $1\le i\le 4$. It follows from the formulas $(v_2\circ \Psi_2)(a_1\ot a_2)=(v_2\circ \Psi_2)(a_2\ot a_3)=0$ and the fact that
$$
(z\circ \Psi_4)(u\ot v\ot a_2\ot a_3)=0
$$
for all $u,v\in A$. We have
$$((z\circ \Psi_4)\circ_i(v_2\circ \Psi_2)\big)(a\ot b\ot y\ot y\ot y)=0$$
 for $i=3$ and $i=4$ because $(v_2\circ \Psi_2)(y\ot y)=v_2(1\ot r_y\ot 1)=0$. Moreover
 $$(z\circ \Psi_4)\big((v_2\circ \Psi_2)(a\ot b)\ot y\ot y\ot y\big)=zt_3\big((v_2\circ \Psi_2)(a\ot b)\Psi_3(y\ot y\ot y)\big)=0.$$
 Also we have
 $$
 \begin{aligned}
&(z\circ \Psi_4)\big(a\ot(v_2\circ \Psi_2)(b\ot y)\ot y\ot y\big)=zt_3\Big(at_2\big(v_2t_1(b\ot y\ot 1)\ot r_y\ot 1\big)\Big)\\
=&\begin{cases}
zt_3(at_2(xyx\ot r_y\ot 1)),&\mbox{if $b=xy$,}\\
zt_3(at_2(xyxy\ot r_y\ot 1)),&\mbox{if $b=yxy$,}\\
0,&\mbox{if $b\in\mathcal{B}\setminus\{xy,yxy\}$}\\
\end{cases}=
\begin{cases}
1,&\mbox{if $a=y$, $b=yxy$,}\\
0,&\mbox{if $a=x$ or $b\in\mathcal{B}\setminus\{yxy\}$.}\\
\end{cases}
\end{aligned}
 $$
 Thus
 $$
 ((z\circ \Psi_4)\circ(v_2\circ \Psi_2)\big)(a\ot b\ot y\ot y\ot y)=\begin{cases}
1,&\mbox{if $a=y$, $b=yxy$,}\\
0,&\mbox{if $a=x$ or $b\in\mathcal{B}\setminus\{yxy\}$.}\\
\end{cases}
 $$

We have
$$((z\circ \Psi_4)\circ_i(v_2\circ \Psi_2)\big)(a\ot b\ot x\ot x\ot x)=(z\circ \Psi_4)(a\ot b\ot x\ot x)$$
 for $i=3$ and $i=4$. So
 $$
 ((z\circ \Psi_4)\circ_3(v_2\circ \Psi_2)\big)(a\ot b\ot x\ot x\ot x)+((z\circ \Psi_4)\circ_4(v_2\circ \Psi_2)\big)(a\ot b\ot x\ot x\ot x)=0.
 $$
 Further we have
 $$
 (z\circ \Psi_4)\big((v_2\circ \Psi_2)(a\ot b)\ot x\ot x\ot x\big)=zt_3(v_2t_1(a{\sf C}(b))\ot 1).
 $$
 Direct calculations show that
 $$
 v_2t_1(a{\sf C}(b))=\begin{cases}
 x,&\mbox{if $a=x$, $b=x$,}\\
  0,&\mbox{if $a=x$, $b\in\{y,yx,yxy\}$ or $a=y$, $b\in\{x,y,xy,xyx,xyxy\}$,}\\
   xy+xyxy,&\mbox{if $a=x$, $b=xy$,}\\
   xyx,&\mbox{if $a=x$, $b=xyx$ or $a=y$, $b=yx$,}\\
   xyxy,&\mbox{if $a=x$, $b=xyxy$ or $a=y$, $b=yxy$.}\\
 \end{cases}
 $$
  Finely we have
 $$
 (z\circ \Psi_4)\big(a\ot (v_2\circ \Psi_2)(b\ot x)\ot x\ot x\big)=zt_3(at_2(v_2t_1(b\ot x\ot 1)\ot r_x\ot 1)).
 $$
Note that $t_3(yt_2(u\ot r_x\ot 1))=0$ for any $u\in \mathcal{B}$.
  Direct calculations show that
 $$
 v_2t_1(b\ot x\ot 1)=\begin{cases}
 x,&\mbox{if $b=x$,}\\
  0,&\mbox{if $b\in\{y,xy,yxy\}$,}\\
   yx+xyxy,&\mbox{if $b=yx$,}\\
    xyx,&\mbox{if $b=xyx$,}\\
     xyxy,&\mbox{if $b=xyxy$.}\\
 \end{cases}
 $$
 Then
$$
\begin{aligned}
&((z\circ \Psi_4)\circ(v_2\circ \Psi_2)\big)(a\ot b\ot x\ot x\ot x)\\
=&\begin{cases}
 1,&\mbox{if $a=x$, $b=xy$ or $a=x$, $b=yx$ or $a=y$, $b=yxy$}\\
 0,&\mbox{if $a=x$, $b\in\mathcal{B}\setminus\{xy,yx\}$ or $a=y$, $b\in\mathcal{B}\setminus\{yxy\}$.}
 \end{cases}
 \end{aligned}
$$
Thus we have
$$
\begin{aligned}
\big((z\circ \Psi_4)\circ(v_2\circ \Psi_2)\big) \Phi_5(1\ot x\ot 1)&=xy+yx,\\
\big((z\circ \Psi_4)\circ(v_2\circ \Psi_2)\big) \Phi_5(1\ot y\ot 1)&=x+x=0.
 \end{aligned}
$$
So $[v_2,z]=(xy+yx,0)=0$, because $(xy+yx,0)\in B^1(A)=B^5(A)$ by Remark \ref{denot_coc_cob} and the 4-periodicity of the resolution $P_*$.
\end{Proof}

We now can prove a theorem which describes the BV structure on $HH^*(A)$.

\begin{Thm}\label{Thm} Let $A=kQ_8$, $\mathrm{char}k=2$ and $\Delta$ be the BV differential from Theorem \ref{Tradler}. Then\\
$1)$ $\Delta$ is equal 0 on the generators of $HH^*(A)$ from Theorem \ref{Gen};\\
$2)$ $\Delta(ab)=0$ for $a\in\{v_1,v_2, v_2'\}$, $b\in\{p_1,p_2,p_2',p_3\}$;\\
$3)$ $\Delta(az)=0$ if $a$ is a generator of $HH^*(A)$ from Theorem \ref{Gen};\\
$4)$ $\Delta$ satisfies the equalities
$$\begin{aligned}
\Delta(p_1u_1)&=\Delta(p_3u_1)=\Delta(p_2u_1')=p_2',\,\Delta(p_2u_1)=\Delta(p_2'u_1')=p_1,\\
\Delta(p_2'u_1)&=\Delta(p_1u_1')=\Delta(p_3u_1')=p_2,\,\Delta(u_1v_1)=\Delta(u_1'v_2')=(u_1')^2+v_2,\\
\Delta(u_1'v_1)&=\Delta(u_1v_2)=u_1^2+v_2',\,\Delta(u_1'v_2)=\Delta(u_1v_2')=v_1.
\end{aligned}$$

Points $1)$--$4)$ with Theorem \ref{Gen} determines BV algebra structure (and in particular Gerstenhabber algebra structure) on $HH^*(A)$.
\end{Thm}
\begin{Proof}
Points 1)--4) follow from Lemmas \ref{BV-structure1}--\ref{BV-structure6}. To determine BV algebra structure we need the value of $\Delta$ on generators and all their pairwise products. Point 1) determines $\Delta$ on generators. Points 2)--4) determine $\Delta$ on all pairwise products of generators except zero products (see Theorem \ref{Gen}) and squares of generators. All the listed products are zero in characteristic two. So BV structure is fully determined.
\end{Proof}

\begin{Cor} Let $A=kQ_8$, $\mathrm{char}k=2$ and $[,]$ be the Gerstehaber bracket from Theorem \ref{Tradler}. Then the bracket is zero for all pairs of generators of $HH^*(A)$ from Theorem \ref{Gen} exept:

$$\begin{aligned}\
[p_1,u_1]&=[p_3,u_1]=[p_2,u_1']=p_2',\,[p_2,u_1]=[p_2',u_1']=p_1,\\
[p_2',u_1]&=[p_1,u_1']=[p_3,u_1']=p_2,\,[u_1,v_1]=[u_1',v_2']=(u_1')^2+v_2,\\
[u_1',v_1]&=[u_1,v_2]=u_1^2+v_2',\,[u_1',v_2]=[u_1,v_2']=v_1.
\end{aligned}$$

This completely determines Gerstehaber algebra structure on $HH^*(A)$.

\end{Cor}

\begin{Proof}
From the Theorem \ref{Thm} we know, that BV-differential $\Delta$ equals zero on any generator of BV-algebra $HH^*(A)$. Then using formula from the definition \ref{Def-BV-algebra} one immediately has $[a,\, b]=\Delta(a\smile b)$ for any $a,\, b$ from the set of generators of $HH^*(A)$.
\end{Proof}

\end{document}